\documentclass[a4paper,12pt]{article} \usepackage{bbm}
\usepackage{amsbsy,amsfonts,amsmath,sectsty,graphicx,cancel,xspace,hyperref,todonotes}

\newcommand{\eg}{{\em e.g.\/}\xspace} 
\newcommand{\ie}{{\em i.e.\/}\xspace}

\newcommand{\fbar}{{\overline F}} 
\newcommand{\pibar}{{\fbar_\Pi}}
\newcommand{\bbl}{\mbox{$\mathbbm 1$}} 
\newcommand{\R}{\mbox{$\mathbbm R$}} 
\newcommand{\E}{\mathbbm E}

\title{Proper scoring rules for survival times}
\author{Philip Dawid}
\date{}

\begin{document}
\maketitle

\abstract{This note presents a general form for a proper scoring rule
  that can be used for possibly right-censored survival data without
  requiring any knowledge or assumptions about the censoring process,
  so long as it is independent of survival time.}


\section{Survival data}
Let $T$ be a possibly censored non-negative survival time.  Censoring
is at a random time $C \leq \infty$, independent of $T$.  We observe
$T$ if $T\leq C$, else we observe $C$ (and so only know $T > C$).
Equivalently, we observe $M = \min\{C,T\}$, and an indicator
$\Delta = \bbl(T \leq C)$, which takes value 1 if $T \leq C$ (\ie, we
have observed $T=M$), or 0 if $T > C$ (in which case we only know that
$T$ was censored at $M$).

\section{A Bregman type scoring rule}
\label{sec:brier}
Let $\psi$ be a concave function on $\R^+$, and define
\begin{equation}
  \label{eq:gamma}
  \gamma(\lambda) := \psi(\lambda) - \lambda\psi'(\lambda).
\end{equation}
(In case $\psi$ is not differentiable at $\lambda$, fix an arbitrary
supergradient to $\psi$ there, and interpret $\psi'(\lambda)$ as the
slope of that supergradient.)

Consider the scoring rule, for when You quote an absolutely
contintuous distribution $Q$ for $T$, and $T=t$:
\begin{equation}
  \label{eq:genpsr}
  S(t,Q) = \int_0^m \gamma\left\{\lambda_Q(u)\right\}\,du + \psi'\left\{\lambda_Q(m)\right\}\,\delta.
\end{equation}
Here $\lambda_Q(u) = q(u)/\fbar_Q(u)$ denotes Your hazard function for
$T$ at $u$, with $q$ the density, and $\fbar_Q = 1 - F_Q$ the survivor
function, of $T$ (all under $Q$).  Also $m$ and $\delta$ are,
respectively, the values of $M$ and $\Delta$ (which -- though the
notation does not make this explicit --- depend on the value $c$ of
$C$, as well as on $t$).

\section{Fixed censoring}
\label{sec:fixed}

We first consider the case that $C$ is a fixed value $c$ (possibly
infinite).  Then when $T=t$ we have the observed values
$m := \min\{t,c\}$ for $M$, and $\delta := \bbl(t\leq c)$ for
$\Delta$.  Suppose You have quoted (absolutely continuous)
distribution $Q$ for $T$; while Your ``true'' distribution is $P$.

In this case, \eqref{eq:genpsr} can be expressed as
\begin{equation}
  \label{eq:genpsrc}
  S_c(t,Q) = \int_0^c  \gamma\left\{\lambda_Q(u)\right\}\bbl(u < t)\,du +\psi'\left\{\lambda_Q(t)\right\}\, \bbl(t\leq c).
\end{equation}
Your expected score, when $T \sim P$, supposed absolutely continuous
with density function $p$, is thus
\begin{equation}
  \label{eq:genspqc}
  S_c(P,Q) = \int_0^c \left[\psi\left\{\lambda_Q(u)\right\} +\left\{\lambda_P(u) - \lambda_Q(u)\right\} \psi'\left\{\lambda_Q(u)\right\} \right]\fbar_P(u)\,du.
\end{equation}
Consequently we have entropy function
\begin{equation}
  \label{eq:genhc} 
  H_c(P) :=S_c(P,P) = \int_0^c \psi\left\{\lambda_P(u)\right\} \fbar_P(u)\,du
\end{equation}
and discrepancy function
\begin{equation}
  \label{eq:gendc} 
d_c(P,Q) := S_c(P,Q)-H_c(P) =  \int_0^c \rho\left\{\lambda_P(u),\lambda_Q(u)\right\} \fbar_P(u)\,du
\end{equation}
where $\rho(a,b) := \psi(b)+(a-b)\psi'(b)-\psi(a) \geq 0$ by concavity
of $\psi$.  It follows that $S$ is a proper scoring rule.

\section{Random censoring}
\label{sec:rand}

Suppose now that $C$ is random, but taken as independent of $T$ under
$P$.  We continue to use scoring rule \eqref{eq:genpsr}.  Then, by
independence, Your overall expected score is
\begin{eqnarray}
  \label{eq:exprand}
  S(P,Q) = \E_{C\sim\Pi} S_C(P,Q)
\end{eqnarray}
where $\Pi$ is Your distribution for $C$.  Since, for any $c$, $S_c$
is minimised in $Q$ at $Q=P$, so too is $S(P,Q)$, irrespective of the
specification of $\Pi$.  Hence \eqref{eq:genpsr} continues to define a
proper scoring rule, even under independent random censoring.  It will
be strictly proper when $\psi$ is strictly concave, and $C$ has
infinite range in the sense that the survivor function $\pibar(u)$ of
$\Pi$ is positive for all finite $u$.

Under random censoring we obtain
\begin{eqnarray}
  \nonumber
  S(P,Q) &=& \int_0^\infty \left[\psi\left\{\lambda_Q(u)\right\} +\left\{\lambda_P(u) - \lambda_Q(u)\right\} \psi'\left\{\lambda_Q(u)\right\} \right]\pibar(u)\fbar_P(u)\,du\\
  \label{eq:genspq}
  &=& \int_0^\infty \left[\psi\left\{\lambda_Q(u)\right\} +\left\{\lambda_P(u) - \lambda_Q(u)\right\} \psi'\left\{\lambda_Q(u)\right\} \right]\fbar_M(u)\,du
\end{eqnarray}
where $\fbar_M$ is the survivor function of $M = \min\{T,C\}$ under
the ``true'' joint distribution of $(T,C)$ (in which we have
independence).  Similarly we have
\begin{eqnarray}
  \label{eq:genh}
  H(P) &=& \int_0^\infty \psi\left\{\lambda_P(u)\right\} \fbar_M(u)\,du\\
  \label{eq:gend}
  d(P,Q) &=& \int_0^\infty \rho\left\{\lambda_P(u),\lambda_Q(u)\right\} \fbar_M(u)\,du.
\end{eqnarray}
Note however that all these quantities depend on the censoring process
currently operating, and might not be relevant for predicting future
behaviour under a different censoring process (or none).

\subsection{Estimation}
\label{sec:est}

For a statistical model $\{P_\theta\}$, we can estimate $\theta$ by
minimising the total empirical score $\sum S(t_i, P_\theta)$.
Typically this will be obtainable by setting to $0$ the derivative
with respect to $\theta$, and general theory shows that this will
yield an unbiased estimating equation --- and typically a consistent
estimator --- irrespective of the (independent) censoring process.
However, the variance of the estimator will depend on that process.

For pre-assigned $Q$, an unbiased estimator of $S(P, Q)$ is,
trivially, the average of $S(t,Q)$ over a random sample of
$(m_i,\delta_i)$ pairs.  This can be used, and is valid, so long only
as we can specify $Q$ --- no knowledge or assumption about $P$ is
required.

On integrating \eqref{eq:genh} by parts\footnote{This is assuming the
  ``boundary terms'' vanish, which will hold so long as
  $\Psi_P(u)\fbar_M(u) \rightarrow 0$ as $u\rightarrow\infty$.}  we
can express it as
\begin{equation}
  \label{eq:inth}
  H(P) = \E\{\Psi_P(M)\}
\end{equation}
where 
\begin{equation}
  \label{eq:Psi}
  \Psi_P(u) := \int_0^u \psi\{\lambda_P(s)\}\,ds.
\end{equation}
Since we observe $M$, we can consistently estimate \eqref{eq:inth} by
the empirical average of the $\Psi_P(m_i)$ --- so long as we can
compute or estimate the function $\Psi_P$, which will typically
require being able to specify the form of $P$, \eg, under an
assumption that $P$ and $Q$ both belong to a common, specified model.
In that case we can also consistently estimate the discrepancy $d(P,Q)
= S(P,Q) - H(P)$ (but again, only for the actual censoring process
underlying the observed data).

\section{Examples}
\label{sec:ex}
\subsection{Brier type scoring rule}

For the case $\psi(x) \equiv -x^2$, we get
\begin{eqnarray*}
  \label{eq:entc}  
  S(t, Q) &=& \int_0^m \lambda_Q(u)^2\,du - 2 \lambda_Q(m)\,\delta\\
  \label{eq:hc}
  S(P, Q) &=& \int_0^\infty  \left\{ \lambda_Q(u)^2 - 2 \lambda_Q(u)\lambda_P(u)\right\}\fbar_M(u)\, du\\
  H(P) &=& - \int_0^\infty  \lambda_P(u)^2\fbar_M(u)\,du\\
  \label{eq:discc}
  d(P,Q) &=&  \int_0^\infty  \left\{\lambda_Q(u)-\lambda_P(u)\right\}^2\fbar_M(u)\,du.
\end{eqnarray*}

\subsection{Tsallis  type scoring rule}

More generally, take $\psi(x) \equiv - x^\beta$ ($\beta> 1$).  We get
\begin{eqnarray*}
  \label{eq:tsentc}  
  S(t, Q) &=& (\beta-1)\int_0^m \lambda_Q(u)^\beta\,du - \beta \lambda_Q(m)^{\beta-1}\,\delta\\
  \label{eq:tshc}
  S(P, Q) &=& \int _0^\infty  \left\{ (\beta-1)\lambda_Q(u)^\beta - \beta \lambda_P(u)\lambda_Q(u)^{\beta-1}\right\}\fbar_M(u)\, du\\
  H(P) &=& - \int_0^\infty  \lambda_P(u)^\beta\,\fbar_M(u)\,du\\
  \label{eq:tsdiscc}
  d(P,Q) &=&  \int_0^\infty  \left\{
    (\beta-1)\lambda_Q(u)^\beta  + \lambda_P(u)^\beta - \beta \lambda_P(u)\lambda_Q(u)^{\beta-1}
  \right\}\fbar_M(u)\,du.
\end{eqnarray*}

\subsection{Log type scoring rule}
\label{sec:log}

For $\psi(x) = x-x\log x$, and assuming $\lambda_Q$ always positive,
we find
\begin{eqnarray}
  \label{eq:psrlog} 
  S(t,Q) =  \Lambda_Q(m) -  \log \lambda_Q(m)\,\delta
\end{eqnarray}
where $\Lambda_Q(m) := \int_0^m \lambda_Q(s)\,ds = -\log \fbar_Q(m)$.

Then we get
\begin{eqnarray*}
  \label{eq:psrcPlog}
  S(P,Q) &=& \int_0^\infty  \left\{\lambda_Q(u)-\lambda_P(u) \log \lambda_Q(u)\right\}\fbar_M(u)\,du\\
  \label{eq:psrchlog}
  H(P) &=&  \int_0^\infty  \left\{\lambda_P(u) - \lambda_P(u)\log \lambda_P(u)\right\} \fbar_M(u)\,du\\
  \label{eq:psrcdlog}
  d(P,Q) &=& \int_0^\infty \left\{\lambda_Q(u)-\lambda_P(u) + \lambda_P(u)\log\frac{\lambda_P(u)}{\lambda_Q(u)}\right\}\fbar_M(u) \,du.
\end{eqnarray*}

\section{Application to exponential distribution}
Suppose we want to fit an exponential distribution $Q_\alpha$, with
density function
\begin{math}
  \alpha e^{-\alpha t},
\end{math}
survivor function
\begin{math}
  e^{-\alpha t},
\end{math}
and constant hazard function $\alpha$.

We minimise the total empirical score:
\begin{eqnarray*}
  \lefteqn{\gamma(\alpha) \sum_im_i + \psi'(\alpha)\sum_i\delta_i}\\
  &=& A\{\psi(r) +\rho(r,\alpha)\}
\end{eqnarray*}
where $A=\sum_im_i$, $r=\sum_i\delta_i/\sum_im_i$.  Thus, for any
choice of concave $\psi$, and any independent censoring process, a
minimum-score estimate (unique if $\psi$ is strictly concave) is
$$\widehat\alpha = r = \frac{\sum_i \delta_i}{\sum_i m_i},$$
which is just the number of failures observed, divided by the total
time on test --- an appealing solution.

\section*{Comment}
The propriety of the scoring rule \eqref{eq:genpsr} was announced
without proof in \S3.5 of \cite{apd/mm:metron}, which may be consulted
for general background.  This version of the note incorporates some
corrections and suggestions made by ChatGPT.

\end{document}